\documentstyle[amscd]{amsart}
\input{bull-art}
\newcommand{\bb}{\Bbb}

\newcommand{\cx}{{\bb C}}

\newcommand{\integers}{{\bb Z}}

\newcommand{\reals}{{\bb R}}

\newcommand{\Bold}[1]{\medskip \noindent {\bf #1 
}\nopagebreak}

\newcommand{\arrow}{\rightarrow}
\newcommand{\bdry}{\partial}
\newcommand{\closure}{\overline}
\newcommand{\compos}{\circ}

\newcommand{\dirsum}{\oplus}

\newcommand{\equi}{\sim}

\newcommand{\mem}{\in}

\newcommand{\superset}{\supset}

\newcommand{\EQ}{\;\;=\;\;}

\newcommand{\st}{\; : \;}         

\newcommand{\chat}{\widehat{\cx}}

\newcommand{\zed}{\integers}

\newcommand{\diam}{\operatorname{diam}}

\newcommand{\interior}{\operatorname{int}}
\renewcommand{\Im}{\operatorname{Im}}

\renewcommand{\mod}{\operatorname{mod}}

\newcommand{\PSL}{\operatorname{PSL}}

\newtheorem{theorem}{Theorem}[section]

\newtheorem{cor}[theorem]{Corollary}

\newcommand{\cP}{{\cal P}}

\newcommand{\hdtwop}{HD2$'$}
\newcommand{\hdtwor}{HD2$\reals$}
\newcommand{\nilftwor}{NILF2$\reals$}

\theoremstyle{remark}

\newtheorem{remarks}{Remarks}

\begin{document}
\def\currentvolume{31}
\def\currentissue{2}
\def\currentyear{1994}
\def\currentmonth{October}
\def\copyrightyear{1994}
\def\currentpages{155-172}

\title{Frontiers in complex dynamics}
 
\author{
        Curtis T. McMullen}
\address{Mathematics Department\\
        University of California\\
        Berkeley, California 94720}
\email{ctm@@math.berkeley.edu}
\thanks{
Based on a lecture
presented to the AMS-CMS-MAA joint meeting,
Vancouver, BC, August 16, 1993.  Supported in
part by the NSF}
\date{February 1, 1994}
\maketitle

\section{Introduction}

Rational maps on the Riemann sphere 
occupy a distinguished niche in the
general theory of smooth dynamical systems.
First, rational maps are complex-analytic,
so a broad spectrum of techniques can contribute
to their study (quasiconformal mappings, potential
theory, algebraic geometry, etc.).
The rational maps of a given degree form a
finite-dimensional manifold, so exploration of
this {\em parameter space} is especially tractable.
Finally, some of the conjectures once proposed 
for {\em smooth} dynamical systems (and now known
to be false) seem to have a definite chance of
holding in the arena of rational maps.

In this article we survey a small constellation
of such conjectures centering around the
density of {\em hyperbolic} rational maps ---
those which are dynamically the best behaved.
We discuss some of the evidence and logic
underlying these conjectures, and sketch
recent progress towards their resolution.

Our presentation entails only a brief
account of the basics of complex dynamics;
a more systematic exposition can be found in
the survey articles 
\cite{Douady:Bourbaki:dynamics},
\cite{Blanchard:survey}, and
\cite{Eremenko:Lyubich:survey};
the recent books
\cite{Beardon:book:rational}
and
\cite{Carleson:Gamelin:book:dynamics};
and Milnor's lecture notes
\cite{Milnor:dynamics:lectures}.

\section{Hyperbolic rational maps}

A {\em rational map} $f : \chat \arrow \chat$ is
a holomorphic dynamical system on the Riemann
sphere $\chat = \cx \cup \{\infty\}$.
Any such map can be written as a quotient
\begin{displaymath}
	f(z) \EQ
	\frac{P(z)}{Q(z)} \EQ
	\frac{
	a_0 z^d + \ldots + a_d 
	}{
	b_0 z^d + \ldots + b_d 
	}
\end{displaymath}
of two relative prime polynomials $P$ and $Q$.
The {\em degree} of $f$ can be
defined topologically or algebraically; 
it is the number of preimages of a typical point $z$,
as well as the maximum of the degrees of $P$ and $Q$.
The fundamental problem in the dynamics of rational
maps is to understand the behavior of high
iterates
\begin{displaymath}
	f^n (z) = 
	\underbrace{
	(f \compos f \compos \ldots \compos f)
		}_{\text{$n$ times}}
	(z) .
\end{displaymath}

Any rational map of degree $d>1$ has both
expanding and contracting features.
For example, $f$ must be expanding on average,
because it maps the Riemann sphere over itself
$d$ times.  Indeed,
with respect to the spherical metric
(normalized to have total area one),
\begin{displaymath}
	\int_{\chat} \|(f^n)'\|^2 
	= d^n \arrow \infty,
\end{displaymath}
so the derivative of $f^n$ is very large on
average.
On the other hand, $f$ has $2d-2$ {\em critical
points} $c$ where $f'(c) = 0$.  Near $c$,
the behavior of $f$ is like that of $z \mapsto z^n$ 
near the origin, for some $n>1$; thus $f$ is
highly contracting near $c$.
Tension between these two aspects of $f$
is responsible for much of the complexity of
rational maps.

To organize these features of $f$, we introduce the
{\em Julia set} $J(f)$ --- the locus of chaotic
dynamics; and the {\em postcritical set}
$P(f)$ --- which contains the ``attractors'' of $f$.
The Julia set can be defined as the closure of the
set of {\em repelling periodic points} for $f$.
Here a point $z$ is {\em periodic} if $f^p(z) =z$ 
for some $p>0$; it is 
\begin{displaymath}
        \begin{array}{ll}
        \text{repelling} & \text{if}\; |(f^p)'(z)|>1,\\
        \text{indifferent} & \text{if}\; |(f^p)'(z)|=1,
        \;\;\text{and}\\
        \text{attracting} & \text{if}\; |(f^p)'(z)|<1.
        \end{array}
\end{displaymath}
The forward orbit $E$ of a periodic point is called
a {\em cycle}, because $f|E$ is a cyclic permutation.

The derivative gives a first approximation to the
behavior of $f^p$ near the periodic point;
for example, all
points in a small neighborhood of an attracting point
$z$ tend towards $z$ under iteration of $f^p$.
On the other hand, a repelling point pushes away
nearby points, so the behavior of forward iterates
is difficult to predict.

The Julia set is also the smallest
closed subset of the sphere such that $|J(f)|>2$ and
$f^{-1}(J) = J$.  Its complement,
$\Omega=\chat-J(f)$, sometimes called the
{\em Fatou set}, is the largest open set such
that the iterates $\langle f^n|\Omega\rangle$ form a
normal family.

The {\em postcritical set} $P(f)$ is the closure of
the forward orbits of the critical points of $f$:
\begin{displaymath}
	P(f) = \closure{ \bigcup_{n>0, \;f'(c) = 0 }
		f^n(c) } .
\end{displaymath}

The postcritical set plays a crucial role with
respect to the attractors of $f$. For example:

\begin{theorem}
\label{thm:attr:cp}
	Every attracting cycle $A$ attracts
a critical point. 
\end{theorem}

\begin{pf}
Let $U = \{z \st d(f^n(z),A) \arrow 0\}$
for the spherical metric; $U$ is open, and
$f^{-1}(U)=U$.  If $U$ contains no critical point,
then $f|U$ is a covering map; but then the
Schwarz lemma implies $f$ is an isometry for
the hyperbolic metric, which is impossible
because $A$ is attracting. 
\end{pf}

Thus $A \subset P(f)$ and the
number of attracting cycles is bounded by
the number of critical points, which in
turn is bounded by $2\deg(f)-2$.

This theorem is of practical as well as theoretical
value.
For example, if $f(z) = z^2 + c$ has an attracting
cycle of period $100$, this cycle can be easily
located as $\lim f^n(0)$; a few million iterates should
yield reasonable accuracy.
(Compare this to the prospect of
computing the $1.27 \times 10^{30}$ roots of the
equation $f^{100}(z) = z$.)

We can now introduce the property of
{\em hyperbolicity}, which will be central in
the remaining discussion.
Let $f$ be a rational map of degree $d>1$.

\begin{theorem}
The following conditions are equivalent\,\rom{:}
\begin{enumerate}
	\item
	All critical points of $f$
	tend to attracting cycles under iteration.
	\item
	The map $f$ is expanding on its Julia set.
That is, there exists a conformal metric 
$\rho$ on the sphere such that $|f'(z)|_\rho > 1$
for all $z \mem J(f)$.
	\item
	The postcritical set and the Julia set are
	disjoint $(P(f)\cap J(f) = \emptyset)$.
\end{enumerate}
\end{theorem}

\Bold{Definition.}
When the above conditions hold, $f$ is
{\em hyperbolic}.  

The Julia set of a hyperbolic rational map is
thin:  its area is zero, and in fact its Hausdorff
dimension is strictly less than two
\cite{Sullivan:CDS}.
Every point outside the Julia set tends towards a
finite attractor $A \subset \chat$:  that is,
the spherical distance $d(f^n(z),A) \arrow 0$
as $n \arrow \infty$.  The set $A$ consists exactly
of the union of the attracting cycles for $f$.

Thus for a hyperbolic rational map, we can predict
the asymptotic behavior of all points in an open,
full-measure subset of the sphere: they converge
to $A$.  

\Bold{Example.}
Figure \ref{fig:hypex} depicts the Julia set of a
rational map of degree $11$.  The Julia set is
in black; its complement contains
twenty large white regions,
ten of which are visible in the picture.
The attractor $A$ consists of one point in the
``center'' of each large white region.
Under iteration,
every point outside the Julia set eventually lands
in one of the large white regions and is then
attracted to its center.
The Julia set is the thin set of 
``indecisive'' points forming the boundary between
regions converging to one point of $A$ or to another.

This rational map is especially symmetric:  it commutes
with the symmetries of the dodecahedron, and it can
be used to solve the quintic equation
(but that is another story; see \cite{Doyle:McMullen}).


We can now state one of the central open problems
in the field.

\begin{figure}[hb]
\vskip18pc
\caption{Julia set of a hyperbolic rational map with
the symmetries of the icoashedron.}\label{fig:hypex}
\end{figure}

\Bold{Conjecture HD.}
{\em
Hyperbolic maps are open and dense among
all rational maps.
}

It is easy to see that hyperbolicity is an open
condition, but the density of hyperbolic dynamics
has so far eluded proof.

Given recent events in number theory, I looked into
the possibility of naming the above conjecture
{\em Fatou's Last Theorem}.  Unfortunately, the name
is unjustified.
Speaking of hyperbolicity, Fatou writes in his
1919--20 memoir \cite[p.73]{Fatou:rational:II}:

\begin{quote}
Il est probable, mais je n'ai pas approfondi la
question, que cette propri\'et\'e appartient \`a
toutes les substitutions g\'en\'erales, c'est-\`a-dire
celles dont les coefficients ne v\'erifient aucune relation
particuli\`ere.\footnote{
I am grateful to Eremenko, Lyubich, and Milnor for
providing this reference.
}
\end{quote}

There is no indication of even a marginal proof.
Moreover, Fatou may have intended by his last statement 
that the non-hyperbolic rational maps should be contained 
in a countable union of proper subvarieties.
This is false, by an elementary argument 
\cite[Proposition 3.4]{Lyubich:stability};
in fact, non-hyperbolic
maps have positive measure among all rational
maps of a given degree \cite{Rees:ergodic:ratl}.

\Bold{Structural stability.}
A pair of rational maps $f$ and $g$ are
{\em topologically conjugate} if there is a
homeomorphism $\phi: \chat \arrow \chat$ such
that $\phi f \phi^{-1} = g$.  A rational map
$f$ is {\em structurally stable} if $f$ is
topologically conjugate to all $g$ in
a neighborhood of $f$.

The following close relative of Conjecture HD is
known to be true:

\begin{theorem}[Ma\~{n}\'e, Sad, Sullivan]
	The set of structurally stable rational
maps is open and dense.
\end{theorem}
\begin{pf*}{Sketch of the proof}
Let $N(f)$ be the number of attracting cycles
of $f$, and let $U_0$ be the set of local maxima
of $N(f)$ in the space of rational maps.
Since attracting cycles persist under small
perturbations and $N(f) \le 2d-2$, the set
$U_0$ is open and dense.  
As $f$ varies in $U_0$, its repelling cycles 
are {\em persistently repelling} --- they cannot
become attracting without increasing $N(f)$.
Tracing the movement of repelling periodic points,
we obtain a topological conjugacy between any
two nearby $f$ and $g$ in $U_0$, defined on a
dense subset of their Julia sets.
By the theory of {\em holomorphic motions},
this map extends continuously to a conjugacy
$\phi : J(f) \arrow J(g)$.

Let $U_1 \subset U_0$ be the set of points
where any critical orbit relations
($f^n(c) = f^m(c')$) are locally constant.
It can be shown that $U_1$ is also open and dense,
and the conjugacy $\phi$ can be extended to
the grand orbits of the critical points over $U_1$.
Finally, general results on holomorphic motions
\cite{Sullivan:Thurston}, \cite{Bers:Royden}
prolong $\phi$ to a conjugacy on the whole sphere. 
\end{pf*}

For details see \cite{Mane:Sad:Sullivan}, 
\cite{McMullen:Sullivan:QCDIII}.

In smooth dynamics, the notion of structural
stability goes back at least to the work of
Andronov and Pontryagin in 1937,
and the problem of the density of structurally
stable systems has been known for some time.
In 1965 Smale showed
that structural stability is {\em not} dense
by giving a counterexample in the space of
diffeomorphisms on a 3-torus \cite{Smale:ss}.
Eventually it was found that neither structural
stability {\em nor} hyperbolicity is dense
in the space of diffeomorphisms, even on
2-dimensional manifolds
(see articles by Abrahams-Smale, Newhouse,
Smale, and Williams in \cite{Chern:Smale:book}).

It is thus remarkable that structural stability
is dense within the space of rational maps;
this fact highlights the special character
of these more rigid dynamical systems.
Given the density of structural stability,
to settle Conjecture HD it suffices to prove
that {\em a structurally stable rational map
is hyperbolic}.

More recent results in smooth
dynamics actually {\em support} Conjecture HD;
the implication
(structural stability) $\implies$ (hyperbolicity)
is now known to hold for $C^1$ diffeomorphisms
\cite{Mane:stability}.

\section{Invariant line fields}

What further evidence can be offered for Conjecture HD?

Theoretical support is provided by a 
more fundamental conjecture, which has its roots in
the quasiconformal deformation theory of rational
maps and relates to Mostow rigidity of
hyperbolic 3-manifolds.  
To describe this conjecture, we will first
give an example of a non-hyperbolic rational map ---
indeed, a rational map whose Julia set is the
entire Riemann sphere.

The construction begins
with a complex torus $X = \cx/\Lambda$, where
$\Lambda = \zed \dirsum \tau \zed$ is a lattice
in the complex plane. 
Choose $n>1$, and
let $F:X \arrow X$ be the degree $n^2$ holomorphic
endomorphism given by $F(x) = nx$.  
Since $|F'(x)| = n>1$, the map $F$ is uniformly
expanding, and it is easy to see that repelling
periodic points of $F$ are dense on the torus $X$.
(For example, all points of order $n^k-1$ in the group
law on $X$ have period $k$ under $F$.)
Thus the Julia set of $F$, appropriately interpreted,
is the whole of $X$.

The quotient of $X$ by the equivalence relation
$x \equi -x$ is the Riemann sphere; the quotient
map $\wp : X \arrow \chat$ can be given by the
Weierstrass $\wp$-function, which presents $X$
as a twofold cover of the sphere branched over
four points.  Since $F(-x) = -F(x)$, the dynamical
system $F$ descends to a rational map 
$f$ such that the diagram
\begin{displaymath}
\begin{CD}
        X & @>{F}>> & X \\
        @V{\wp}VV &     & @V{\wp}VV\\
        \chat & @>{f}>> & \chat
\end{CD}
\end{displaymath}
commutes.

The mapping $f$ can be thought of as an analogue
of the multiple angle formulas for sine and
cosine, since $f(\wp(z)) = \wp(nz)$. 

\Bold{Definition.}
A rational map $f$ is {\em covered by an integral
torus endomorphism} if it arises by the construction
above.\footnote{
This construction goes back to
Latt\`es \cite{Lattes:example}.}

Here are some remarkable features of $f$:
\par
(1)
The Julia set $J(f) = \chat$.  This follows 
easily from the density of repelling periodic
points for $F$ on $X$.
\par
(2)
The mapping $f$ is not rigid:  that is, by
deforming the lattice $\Lambda$ (varying $\tau$),
we obtain a family of rational maps which are
topologically conjugate but not conformally
conjugate.
\par
(3)
Most importantly, the Julia set $J(f)$ carries an
{\em invariant line field}.  
To visualize this line field, first note that
the map $z \mapsto n z$ 
preserves the family of horizontal lines in the 
complex plane.  Thus $F$ preserves the images of such
lines on the torus.  The quotient line family 
turns out to be a foliation by
parallel simple closed geodesics (with respect to
the obvious Euclidean metric) on the torus.
Finally, $f$ preserves the image of this foliation
on the sphere.

Of course there is no way to comb a sphere, so
the image foliation has singularities: 
there are four singular points at the
four critical values of $\wp$.

More generally, an
{\em invariant line field} for $f$, defined on 
a measurable set $E \subset \chat$, is the choice of
a one-dimensional real subspace $L_z$ 
in the tangent space $T_z \chat$ for all
$z \mem E$, such that:
\begin{enumerate}
	\item
$E$ has positive area,
	\item
$f^{-1}(E) = E$,
	\item
the slope of $L_z$ varies measurably with respect to
$z$, and
	\item
the derivative $f'$ transforms $L_z$ into 
$L_{f(z)}$ for all $z$ in $E$.
\end{enumerate}
If $E \subset J(f)$, we say $f$ admits an invariant
line field {\em on its Julia set}.  
Thus the Julia set must have positive measure
before it can carry an invariant line field.

\Bold{Conjecture NILF.}
{\em
A rational map $f$ 
carries no invariant line field on its Julia set,
except when $f$
is covered by an integral torus endomorphism.
}

This conjecture is stronger than the density
of hyperbolic dynamics:

\begin{theorem}[Ma\~n\'e, Sad, Sullivan]
NILF $\implies$ HD.	
\end{theorem}

See \cite{Mane:Sad:Sullivan}, 
\cite{McMullen:Sullivan:QCDIII}.

One attractive feature of conjecture NILF 
is that it shifts the focus of study from 
the family of {\em all}
rational maps to the ergodic theory
of a {\em single} rational map.

In support of this conjecture, and hence of
the density of hyperbolic dynamics, we state
a parallel result for
degree one rational maps.  
Of course a single degree 
one rational map is not very complicated.
The degree one 
mappings form a group, isomorphic to
$\PSL_2(\cx)$, and
the group structure makes it easy to iterate
a single mapping.  To make the dynamical system
more interesting, let us consider more
generally finitely generated
{\em subgroups} $\Gamma \subset \PSL_2(\cx)$
and define the
Julia set $J(\Gamma)$ as the minimal closed
invariant set with $|J|>2$.
We then have:

\begin{theorem}[Sullivan]
\label{thm:nilf}
A discrete finitely generated group 
$\Gamma$ of degree
one rational maps carries no invariant
line field on its Julia set.
\end{theorem}

This result is a (thinly disguised) version
of Mostow rigidity for 
hyperbolic 3-manifolds and orbifolds
with finitely generated fundamental group
\cite{Sullivan:linefield}; in more traditional
terminology, $\Gamma$ is a Kleinian group and
$J(\Gamma)$ is its limit set.

If we allow
$\Gamma$ to be an {\em indiscrete} subgroup of
$\PSL_2(\cx)$, the theorem fails, but in a completely
understood way. For example, the group
\begin{displaymath}
	\Gamma \EQ
	\langle
z\mapsto z+1, z\mapsto z+\tau,z\mapsto n z \rangle,
\end{displaymath}
with $\Im(\tau) > 0$ and $n>1$, has $J(\Gamma) = \chat$, and
it preserves the field of horizontal lines in $\cx$.
This example is simply the universal cover of a torus
endomorphism; in a sense, the exceptions proposed 
in conjecture NILF correspond to the (easily classified)
case of indiscrete groups.

With this result to guide us, why has the 
no-invariant-line-field conjecture remained elusive?
The main reason is perhaps that all
rational maps of degree one lie in a finite-dimensional
Lie group.  This group provides a good geometric
portrait of an arbitrary degree one
transformation.  By contrast, the degree of a
general rational map can tend to infinity
under iteration,
and it is much more difficult to visualize and
control the behavior of a rational map of high degree.

\section{Quadratic polynomials}

The simplest rational maps, apart from those of
degree one, are the {\em quadratic polynomials.}
To try to gain insight into the general theory 
of rational maps, much effort has been devoted
to this special case.  The quadratic
polynomials are remarkably rich in structure, and
many fundamental difficulties are
already present in this family.

From the point of view of dynamics, every
quadratic polynomial occurs exactly once in
the family
\begin{displaymath}
	f_c(z) = z^2 + c\qquad (c \mem \cx),
\end{displaymath}
so the quadratic parameter space can be
identified with the complex plane.

Restricting attention from rational maps to
quadratic polynomials,
it is natural to formulate the following
conjectures.

\Bold{Conjecture HD2.}
{\em 
Hyperbolic maps are dense
among quadratic polynomials.}

\Bold{Conjecture NILF2.}
{\em A quadratic polynomial admits no
invariant line field on its Julia set.}

It turns out that these two conjectures are
{\em equivalent}.  This equivalence is further
evidence for the fundamental nature of
the question of invariant line fields.

Note that $f_c$ has only one critical point in
the complex plane, namely, $z=0$.  Consequently:

\begin{theorem}
	The map $f_c(z) = z^2 + c$ is hyperbolic if and only
if $f^n_c(0) \arrow \infty$ or $f_c$ has an attracting 
periodic
cycle in the finite plane $\cx$.
\end{theorem}

This theorem motivates the following:

\Bold{Definition.}
The {\em Mandelbrot set} $M \subset \cx$ is the set
of $c$ such that $f^n_c(0)$ stays bounded
as $n \arrow \infty$.


\begin{figure}
\vskip14.5pc
\caption{The boundary of the Mandelbrot set.}
\label{fig:M}
\end{figure}

The Mandelbrot set is compact, connected, and
full (this means $\cx-M$ is also connected).
The interior of $M$ consists of countably
many components, appearing as bounded white regions in
Figure \ref{fig:M}.
Thus $M$ can be thought of as a ``tree with
fruit'', the fruit being the components of
its interior (cf. \cite{Douady:descriptions}).

Where does the fruit come from?

\Bold{Conjecture \hdtwop.}
{\em If $c$ lies in the interior of the
Mandelbrot set, then $f_c(z)$ has an
attracting cycle.}

It turns out that hyperbolicity is infectious --- if $U$ is
a component of the interior of the Mandelbrot set, and
$f_c$ is hyperbolic for one $c \mem U$, then $f_c$ is
hyperbolic for {\em all} $c \mem U$.  In this case we
say $U$ is a {\em hyperbolic component} of $\interior(M)$.

It follows that Conjecture
\hdtwop ~is also equivalent to Conjecture
HD2,  so several natural
conjectures concur in the setting of complex quadratic
polynomials.

\Bold{Real quadratic polynomials.}
The Mandelbrot set meets in the real axis
(the horizontal line of symmetry in Figure
\ref{fig:M}) in the interval $[-2,1/4]$.
We can further specialize the conjectures 
above to real quadratics, obtaining:

\Bold{Conjecture \hdtwor.}
{\em Hyperbolicity is dense among real quadratic
polynomials.}

\Bold{Conjecture \nilftwor.}
{\em A real quadratic polynomial admits no invariant
line field on its Julia set.}

The real quadratic polynomials are of special interest
for several reasons.  
First, there are many dynamical
systems which can be roughly modeled on such a
polynomial:  the economy, animal populations, 
college enrollment, etc.
To explain this,
it is convenient to conjugate a
real polynomial in $M$ to the form
$g(x) = \lambda x (1-x)$, where
$0 < \lambda \le 4$, so $g : [0,1] \arrow [0,1]$.
Then $g$ is a ``unimodal map'':  for $x>0$ small
$g(x)$ grows as $x$ grows; but after $x$
passes a critical point ($x=1/2$),
$g$ decreases as $x$ increases.  Thus $g$ might
describe the boom and bust of economic cycles or
the behavior of population from one year to the
next when faced with limited resources.
See Figure \ref{fig:unimodal}, which plots $g$
together with the diagonal $y=x$ and shows an
example where the critical point has period 8.

One can imagine that real numbers correspond to
real life, and one goal of the complex theory is
ultimately to contribute to the understanding of
dynamics over the reals.

\begin{figure}
\vskip12.5pc
\caption{Real quadratic with periodic critical point.}
\label{fig:unimodal}
\end{figure}

Second, some of the combinatorial and geometric
analysis of a quadratic polynomial becomes especially
tractable over the real numbers, because of the
order structure on the real line.  For example,
the forward orbit $\langle f^n_c(0)\rangle$ of the 
critical point
is real when $c \mem \reals$, so the postcritical
set $P(f_c)$ is thin and cannot double back on
itself.

Finally, if we consider $z^2+c$ with both $z$ and
$c$ real,
we can conveniently draw two-dimensional pictures
displaying dynamical features on the $z$ line
as the parameter $c$ varies.  

One such classic computer experiment is the following.
For $c \mem \reals$, let
\begin{displaymath}
	A_c \EQ
	\text{\{limit points of $f^n_c(0)$ as
	$n \arrow \infty$\}}
	\subset \reals
\end{displaymath}
denote the ``attractor'' of $f_c$.  If $f_c$ has an
attracting cycle, then $A_c$ will be equal to that 
finite set.
On the other hand, if $A_c$ is an infinite, then
$f_c$ cannot be hyperbolic.

Now draw the set
$\{(x,c) \st x \mem A_c\}$ as $c$ varies along the
real axis in the Mandelbrot set in the
negative direction, starting just within
the main cardioid; the result appears in 
Figure \ref{fig:cascade}.  
In the main cardioid of $M$, $f_c$ has an
attracting

\begin{figure}[hb]
\vskip17pc
\caption{Bifurcation diagram.}
\label{fig:cascade}
\end{figure}

\noindent
fixed point, so at the bottom of the
figure $A_c$ consists of a single point.
As $c$ decreases, this point bifurcates to a
cycle of period 2, which in turn bifurcates
to period 4, 8, and so on.  Above this ``cascade
of period doublings'' the structure becomes very
complicated and the picture is much darker; there
are large sets of $c$ such that $A_c$ is an entire
interval and the corresponding map is far from
hyperbolic.  The top of the figure corresponds to
$c=-2$.

A blowup of the region near $c=-2$ appears
in Figure \ref{fig:chaos}.  (The prominent smooth
curves in this picture come from the 
forward orbit of the critical point.)
This picture makes apparent the ubiquity of
chaotic dynamics.  In fact, we have the following
result:

\begin{figure}
\vskip17pc
\caption{Blowup near $c=-2$.}
\label{fig:chaos}
\end{figure}

\begin{theorem}[Jakobson]
The set of non-hyperbolic maps has positive measure
in the space of real quadratic polynomials.
\end{theorem}

See \cite{Jakobson:chaos}, \cite{Yoccoz:Bourbaki:Henon}.
Jakobson also shows
that $c=-2$ is a one-sided point of density of
the set of non-hyperbolic maps.

On the other hand, some narrow horizontal windows 
of white are also visible in Figure \ref{fig:chaos};
these ``eyes in the storm'' correspond to hyperbolic
maps, and successive blowups support Conjecture 
\hdtwor:  the hyperbolic windows are apparently dense.

The coexistence of these phenomena leads us to
propose the following:

\Bold{Challenge Question.}
{\em Does $f(z) = z^2-1.99999$ have an attracting
periodic point}?

It is unlikely this question will ever be rigorously
settled, for by Jakobson's theorem, the answer is
almost certainly ``no''.
On the other hand, if
hyperbolicity is indeed dense among real quadratics, then 
we can
change the constant $1.99999$ somewhere past its
trillionth decimal place to obtain a new conjecture
where the answer is ``yes''.
It is hard to imagine a proof that would distinguish
between these two cases.

We can sum up the conjectures put forth so far,
and known implications between them, 
in the following table.  

\begin{center}
\begin{tabular}{|l||c|c|}\hline
~ & Hyperbolic & No invariant \\ 
~ & maps are dense & line fields \\ \hline
Rational maps & HD $\impliedby$ & NILF \\ 
~ & ~ & $\Downarrow$ \\
Quadratic polynomials & HD2 $\iff$ & NILF2 \\ 
~ & ~ & $\Downarrow$ \\
Real quadratic polynomials & \hdtwor $\implies$ & 
\nilftwor \\ \hline
\end{tabular}
\end{center}

\bigskip
Remarkably, the fundamental conjectures concerning
quadratic polynomials (real or complex) can be
subsumed into the following topological statement:

\Bold{Conjecture MLC.}
{\em The Mandelbrot set is locally connected.}

\begin{theorem}[Douady-Hubbard]
\label{thm:DH}
MLC $\implies$ HD\rom{2}, 
HD\rom{2}\<$\Bbb{R}$, NILF\rom{2}, and%
NILF\rom{2}\<$\Bbb R$.
\end{theorem}

Why is locally connectivity such a powerful
property?
One answer comes from a theorem of
Carath\'eodory, which states that the
Riemann mapping 
\begin{displaymath}
\psi  : (\cx-\closure{\Delta}) \arrow (\cx-M)
\end{displaymath}
extends to a continuous map $S^1 \arrow \bdry M$ if
and only if $\bdry M$ is locally connected.
(Here $\Delta$ is the unit disk and $S^1 = \bdry \Delta$.
The Riemann mapping is normalized so that
$\psi(z)/z \arrow 1$ as $z \arrow \infty$.)

If $M$ is locally connected, then each point
$\exp(2\pi i t) \mem S^1$ corresponds to a
unique point $c$ in $\bdry M$.
The {\em external angle} $t$ is a sort of generalized
rotation number, and indeed the mappings corresponding
to rational values of $t$ are well understood.
On the other hand,
the combinatorics of $f_c$ determines the
(one or more) external angles $t$ to which
it corresponds.

If $M$ is locally connected,  then a
quadratic polynomial $f_c$ with $c \mem \bdry M$
is {\em determined} by its
combinatorics, even for irrational external angles.
Using this information,
one can build an abstract 
model for $M$ which is topologically correct; since the
density of hyperbolicity is a topological notion,
it suffices to check it in the abstract model, and
Conjectures HD2 and \hdtwor ~follow.

\section{Renormalization}

We next present some recent breakthroughs in
the direction of the conjectures above.
To explain these results, we will need the concept
of {\em renormalization}.

The local behavior of a rational map can sometimes
be given a linear model.
For example, near a repelling fixed point $p$
with $f'(p) = \lambda$,
one can choose a complex coordinate $z$
so that the dynamics take the form 
$f: z \mapsto \lambda z$.

Renormalization is simply
{\em nonlinear linearization}\,;
that is, one looks for a local model of the dynamics which 
is
a polynomial of degree {\em greater} than one.
We will make this precise in the context of
quadratic polynomials.

Let $f(z) = z^2+c$ with $c$ in the Mandelbrot set.
An iterate $f^n$ is {\em renormalizable} if there
exist disks $U$ and $V$ containing the origin,
with $\closure{U}$ a compact subset of $V$,
such that (a) $f^n : U \arrow V$ is a proper map
of degree two and (b) $f^{nk}(0) \mem U$ for all
$k>0$.
This means that although $f^n$ is a polynomial of
degree $2^n$, it behaves like a polynomial of degree
two on a suitable neighborhood of the critical
point $z=0$.
The restriction $f^n : U \arrow V$ is called a
{\em quadratic-like map}.
A fundamental theorem 
of Douady and Hubbard asserts that any
quadratic-like map
is topologically conjugate to a quadratic polynomial
$g(z) = z^2+c'$;
condition (b) implies $c'$ lies
in the Mandelbrot set and, 
with suitable normalizations, is unique
\cite{Douady:Hubbard:polylike}.

The concept of renormalization explains much of the
self-similarity in the Mandelbrot set and in the
bifurcation diagram for real quadratic polynomials.
For example, there is a prominent window of white
in the midst of the chaotic regime of Figure
\ref{fig:cascade}; a blow-up of this region
appears in Figure \ref{fig:win}.  Remarkably,
three small copies of the entire bifurcation diagram
appear.  This is explained by the fact that
$f^3$ is renormalizable for all values of $c$ in
this window.  As $c$ traverses the window, the
quadratic-like maps $f^3_c : U_c \arrow V_c$
recapitulate the full family of bifurcations of
a quadratic polynomial.
(In the Mandelbrot set, one finds a small homeomorphic
copy of $M$ framing this window on
the real axis.)

\begin{figure}
\vskip17pc
\caption{Recapitulation of bifurcation.}
\label{fig:win}
\end{figure}

\Bold{Infinite renormalization.}
A quadratic polynomial $f$ is {\em infinitely
renormalizable} if $f^n$ is renormalizable for
infinitely many $n>1$.

The prime example of an infinitely renormalizable
mapping is the {\em Feigenbaum polynomial} 
$f(z) = z^2-1.401155\ldots$.  For this map,
a suitable restriction of 
$f^2$ is a quadratic-like map topologically conjugate
to $f$ itself.  It follows that $f^{2^n}$ is
renormalizable for every $n\ge 1$.
Its attractor $A_c$ is a Cantor set representing 
the limit of
the cascade of period doublings visible in
Figure \ref{fig:cascade}.
This Cantor set, the map $f$, and the cascade of
period doublings all exhibit remarkable 
universal scaling features that physicists associate
with phase transitions and that have studied for
many years (see, e.g. the collection 
\cite{Cvitanovic:book:chaos}).

Techniques from complex analysis and Teichm\"uller
theory have been brought to bear by Sullivan 
to provide a conceptual understanding of
this universality \cite{Sullivan:renormalization}.
At the moment the theory applies
only to {\em real quadratics}, that is,
$z^2+c$ with $c \mem \reals$; however, there is little
doubt that universality exists over
$\cx$ \cite{Milnor:hairiness}.

Infinitely renormalizable mappings are
very special.  Remarkably, great progress
has been made towards understanding all other
quadratic polynomials and settling for them 
the conjectures discussed in this paper.
The central result is:

\begin{theorem}[Yoccoz]
\label{thm:mlc}
If $c$ belongs to the Mandelbrot set, then either\,\RM:
\begin{enumerate}
\item[]	$f_c(z) = z^2+c$ is infinitely renormalizable,
	or
\item[] $J(f_c)$ admits no invariant line field and 
	$M$ is locally connected at $c$.
\end{enumerate}
\end{theorem}

Yoccoz's theorem was anticipated by a
breakthrough in cubic polynomials due to Branner
and Hubbard \cite{Branner:Hubbard:cubicsII},
and we will use their
language of tableaux to describe Yoccoz's proof.
(See also
\cite{Milnor:local:connectivity},
\cite{Hubbard:local:connectivity}, and
\cite{Yoccoz:local:connectivity}.)

\begin{pf*}{Sketch of the proof}
Suppose $c \mem M$.
Let  $K(f_c)$ denote the {\em filled Julia set},
that is, the set of $z \mem \cx$ which
remains bounded under iteration of $f_c$; its boundary is 
the
Julia set, and $K(f_c)$ is connected.

Let 
\begin{displaymath}
	\phi_c : (\cx-\closure{\Delta}) \arrow (\cx-K(f_c))
\end{displaymath}
be the Riemann mapping, normalized so that $\phi_c'(z)=1$
at infinity.
It is easy to see that 
\begin{displaymath}
	\phi_c(z^2) = f_c(\phi_c(z));
\end{displaymath}
in other words, $\phi_c$ conjugates the 
$z^2$ to $f_c$.

An {\em external ray} $R_t$ is the image of the
ray $(1,\infty) \exp(2\pi i t)$ under the mapping $\phi_c$;
similarly, an {\em external circle} $C_r$ 
(also called an equipotential)
is the image of $\{z \st |z| = r \}$.
Note that $f_c(R_t) = R_{2t}$ and $f_c(C_r) = C_{r^2}$ by
the functional equation for $\phi_c$.

The main case of the proof arises when all
periodic cycles of $f_c$ are repelling;
let us assume this.
The first step is to try to show that the {\em Julia set}
$J(f_c)$ is locally connected.
To this end,
Yoccoz constructs a sequence $\langle \cP_d\rangle$ of
successively finer tilings of neighborhoods of 
$J(f_c)$.  

To illustrate the method, consider the special case 
$c=i$.  
For this map, the external rays
$R_{1/7}$, $R_{2/7}$, and $R_{4/7}$ converge to a repelling
fixed point $\alpha$ of $f_c$.  
These rays cut the disk bounded by
the external circle $C_2$ into
three tiles (see Figure \ref{fig:puzzle}) called the
{\em puzzle pieces} $\cP_0$ at level $0$.  The pieces
at level $d+1$ are defined inductively as the components
of the preimages of the pieces $\cP_d$ at level $d$.
The new pieces fit neatly inside those already defined,
because the external rays converging to $\alpha$ are
forward-invariant.

The puzzle pieces provide connected neighborhoods of points
in the Julia set.  To show $J(f_c)$ is locally connected,
it suffices to show that $\diam(P_d) \arrow 0$ for any 
nested
sequence of pieces $P_0 \superset P_1 \superset P_2\supset 
\ldots$.

Now
$\diam(P_i) \arrow 0$ will follow if
we can establish 
\begin{displaymath}
	\sum \mod(P_i-P_{i+1}) = \infty;
\end{displaymath}
here each region $P_i-P_{i+1}$ is a
(possibly degenerate) annulus,
and the {\em modulus} $\mod(A) = m$
if the annulus $A$ is 
conformally isomorphic to the standard round annulus
$\{z \st 1 < |z| < \exp(2\pi m)\}$.

\begin{figure}
\vskip16pc
\caption{The Yoccoz puzzle.}\label{fig:puzzle}
\end{figure}

The modulus is especially useful in holomorphic
dynamics because it is invariant under conformal mappings;
more generally, $\mod(A') = \mod(A)/d$ if $A'$ is a
$d$-fold covering of an annulus $A$.

Since the image of a puzzle piece of depth $d>0$ under $f_c$
is again a puzzle piece, the moduli of the
various annuli that can be
formed satisfy many relations.  Roughly speaking, the 
tableau
method allows one to organize these relations and test for
divergence of sums of moduli.

For degree two polynomials, the method succeeds {\em unless}
certain annuli are repeatedly covered by degree two.
Unfortunately, this exceptional case
leads to the convergent sum $1+1/2+1/4+\ldots$,
and so it does not prove local connectivity.
However, one finds this case {\em only} occurs when the 
polynomial
is renormalizable.  

The case of a finitely renormalizable map
$f_c$
can be handled by respecifying the initial tiling
$\cP_0$.  Thus the method establishes locally
connectivity of $J(f_c)$ {\em unless} the mapping
is infinitely renormalizable.

It is a metatheorem that the structure of the
Mandelbrot set at $c$ reflects properties of the 
Julia set $J(f_c)$.
In this case the proof of locally connectivity of
$J(f_c)$ can be adapted, with some difficulty,
to establish locally connectivity of $M$ at
$c$.
A variant of Theorem \ref{thm:DH} then shows 
$J(f_c)$ admits no invariant line field as an
added bonus. 
\end{pf*}

Our own work addresses the infinitely
renormalizable case.
The main result of \cite{McMullen:real} is:

\begin{theorem}
\label{thm:real}
	If $f(z) = z^2+c$ is an infinitely
renormalizable {\em real} quadratic polynomial,
then $J(f)$ carries no invariant line field.
\end{theorem}

When combined with Yoccoz's result, this theorem
implies a positive resolution to Conjecture
\nilftwor, which we restate as follows:

\begin{cor}
	Every component $U$ of the interior of
	the Mandelbrot set that meets the real
	axis is hyperbolic.
\end{cor}

In other words, if one runs the real axis 
through $M$, then all the fruit which is skewered
is good.

\begin{pf*}{Sketch of the proof of Theorem \ref{thm:real}}
By techniques of Sullivan 
\cite{Sullivan:renormalization}, 
the postcritical set 
$P(f)$ of an infinitely renormalizable
real quadratic polynomial is a Cantor set with
gaps of definite proportion at infinitely
many scales.  Using this information and
{\em abandoning} the notion of a quadratic-like
map, we construct instead infinitely many proper
degree two maps $f^n : X_n \arrow Y_n$
(where we do {\em not} require that $X_n \subset Y_n$.)
These maps range in a compact
family up to rescaling.  By the Lebesgue
density theorem, any measurable line field $L_z$
looks nearly parallel on a small enough scale;
using the dynamics, we transport this nearly 
parallel structure to $Y_n$ and pass to the limit.
The result is a mapping with a critical point which
nevertheless preserves a family of parallel lines,
a contradiction.  Thus the original map carries no
invariant line field on its Julia set. 
\end{pf*}

\begin{remarks}
In part, the structure of the argument parallels 
Sullivan's proof of Theorem \ref{thm:nilf};
compactness of the mappings $f^n : X_n \arrow Y_n$
is a replacement for the 
finite-dimensionality of the group of M\"obius
transformations.

The proof also applies to certain
complex quadratic polynomials, those which we call
{\em robust}.  For these maps, the notion of 
``definite gaps'' in the postcritical Cantor set
is replaced by a condition on the hyperbolic lengths of
certain simple closed curves on the Riemann
surface $\chat-P(f)$.

Unfortunately, it is likely that robustness 
can fail for $z^2+c$ when $c$ is allowed to be complex.
Counterexamples can probably be found using a
construction of Douady and Hubbard, which also produces
infinitely renormalizable
quadratic polynomials whose Julia sets are
{\em not} locally connected
\cite{Milnor:local:connectivity}.
\end{remarks}

\section{Further developments}

To conclude, we mention three of the many other
recent developments in complex dynamics
which are most closely connected to the present
discussion.

First,  \'{S}wi\c{a}tek has announced a proof of
Conjecture \hdtwor, the density of hyperbolic
maps in the real quadratic family
\cite{Swiatek:dense}.
This remarkable result settles the topological
structure of bifurcations of real quadratic polynomials.
Note that Conjecture \hdtwor ~implies 
Theorem \ref{thm:real}.

Second, Lyubich has announced a proof of the
local connectivity of the Mandelbrot set at
a large class of infinitely renormalizable points
\cite{Lyubich:lc}.
Thus it seems likely that Conjecture MLC itself is
not too far out of reach.
This conjecture, once settled, will complete our
topological picture of the space of {\em complex}
quadratic polynomials.

\begin{figure}
\vskip15.5pc
\caption{Log of the Mandelbrot set.}
\label{fig:log}
\end{figure}

Finally, Shishikura has settled a long-standing
problem about the {\em geometry} of
the Mandelbrot set by proving that 
$\bdry M$
(although it is topologically one-dimensional)
has Hausdorff dimension two \cite{Shishikura:dim:two}.
To illustrate the complexity of the boundary of
the Mandelbrot set,
Figure \ref{fig:log} renders the image of 
$\bdry M$ under
the transformation $\log(z-c)$ for a certain 
$c \mem \bdry M$.\footnote{
Namely, $c=-0.39054087\ldots -0.58678790i\ldots$,
the point on the boundary of the main cardioid
corresponding to the golden mean Siegel disk.}
Note the cusp on the main cardioid in the upper
right;
looking to the left in the figure corresponds
to zooming in towards the point $c$. 
(It is unknown at this time if $\bdry M$ has
positive area; although the figure looks quite
black in some regions, upon magnification these 
features resolve into fine filaments, apparently of
area zero. Cf. \cite{Milnor:hairiness}.)

In spite of these results and many others,
the main conjectures in complex dynamics
remain open.
Our understanding of parameter space
decreases precipitously
beyond the setting of quadratic polynomials,
and the realm of general rational maps
contains much uncharted territory.
For approaches to cubic polynomials and degree
two rational maps, see
\cite{Milnor:cubics},
\cite{Milnor:quadratic:rational},
\cite{Rees:components},
\cite{Rees:degtwo:I},
\cite{Branner:Hubbard:cubicsI},
and
\cite{Branner:Hubbard:cubicsII}.



\end{document}